\documentclass[11pt,twoside]{amsart}

\usepackage{amssymb, palatino}
\usepackage{color}
\usepackage{ifthen}
\usepackage{euscript}

\date{\today}

\newtheorem{dfn}{Definition}[section]
\newtheorem{thm}{Theorem}
\newenvironment{MainTheorem}
{\trivlist
\item[]\noindent{\bf Main Theorem.}\enspace\ignorespaces}
{\endtrivlist}

\newtheorem{lem}[dfn]{Lemma}
\newtheorem{crl}[dfn]{Corollary}
\newtheorem{prop}[dfn]{Proposition}

\def\inv{^{-1}}

\def\S{\ensuremath{\mathcal S}}
\def\O{\ensuremath{\mathcal O}}

\def\forall{\text{for all }}

\def\P{\ensuremath{\mathbb P}}
\def\N{\ensuremath{\mathbb N}}

\def\R{\ensuremath{\mathbb R}}
\def\Z{\ensuremath{\mathbb Z}}
\def\H{\ensuremath{\mathbb H}}
\def\C{\ensuremath{\mathbb C}}
\def\T{\ensuremath{\mathbb T}}

\def\CP1{\ensuremath{\mathbb C \mathbb P^1}}
\def\Cstar{\ensuremath{\C^*}}
\def\Cs{\ensuremath{\C^*}}
\def\Cn{\ensuremath{\C^n}}

\def\Pn-1{\ensuremath{\P^{n-1}}}

\def\cf{{\em cf. }}

\newcommand{\norme}[1]{|\! | #1 |\! |}

\newcommand{\applic}[5]{
  \ensuremath{\begin{array}[t]{c@{}c@{}c@{}l}
           #1:\ & #2 & \rightarrow & #3 \\
               & #4 & \mapsto     & #5
         \end{array}
        }
}

\newcommand{\innerp}[2]{
\pmb{<} #1 , #2 \pmb{>}
                        }

\newcommand{\setofst}[2]
{ \ensuremath{ \big\{\, #1\ \big\vert\ #2 \big\}
            }
}

\newcommand{\setofstquad}[2]
{ \ensuremath{ \big\{\, #1\ \big\vert\quad #2 \big\}
            }
}

\newcommand{\vect}[1]{ {\scriptstyle{
\left[
\begin{array}{c}
#1
\end{array}
\right]               }}}

\definecolor{red}{rgb}{.7,0,0}
\definecolor{green}{rgb}{0,.7,0}

\newboolean{old_draft} \setboolean{old_draft}{false}
\newcommand{\cmt}[1]
{\ifthenelse{\boolean{old_draft}}{\marginpar{{\sc \scriptscriptsize \color{green} #1}}}{}}
\newcommand{\inred}[1]
{\ifthenelse{\boolean{old_draft}}{{\color{red} #1}}{#1}}

\newboolean{draft}

\newcommand{\margincmt}[1]
{\ifthenelse{\boolean{draft}}{\marginpar{{\sc \scriptsize \color{green} #1}}}{}}
\newcommand{\modifinred}[1]
{\ifthenelse{\boolean{draft}}{{\color{red} #1}}{#1}}

\def\k{\textbf{\large \emph{k}}}
\def\sk{\textbf{\emph{k}}}
\def\z{\textbf{\large \emph{z}}}
\def\sz{\textbf{\emph{z}}}

\def\ub{\textbf{\large \emph{u}}}
\def\vb{\textbf{\large \emph{v}}}
\def\kd{\ensuremath{\tilde{ \textbf{\large \emph{k}} }}}
\def\skd{\ensuremath{\tilde{ \textbf{\emph{k}} }}}
\def\zd{\ensuremath{\tilde{ \textbf{\large \emph{z}} }}}

\def\bo{\textbf{\large \emph{b}}_{\mathbf{1}}}
\def\bd{\textbf{\large \emph{b}}_{\textbf{\emph{d}}}}
\def\bi{\textbf{\large \emph{b}}_{\textbf{\emph{i}}}}

\def\Csd{\ensuremath{(\Cs)^d}}
\def\S{\ensuremath{\mathcal{S}}}

\def\cusp{\ensuremath{\mathcal{C}}}

\def\gldz{\ensuremath{GL_d(\Z)}}

\def\gldzrbo{\ensuremath{GL_d^{\rho > 1}(\Z)}}

\def\EMm{\ensuremath{E_m}}
\def\wt{\ensuremath{\tilde{ w }}}

\title{Holomorphic functions on bundles over annuli}
\author{Dan Zaffran}
\thanks{zaffran@fudan.edu.cn\\
\hspace*{3.5mm} Fudan University, Shanghai.\\
\hspace*{3.5mm} Academia Sinica, Taipei.}

\setboolean{draft}{false}

\begin{document}

\margincmt{LIST OF MODIFICATIONS: 
(page and line numbers with respect to *previous* 
version are indicated in brackets)\\ 
P.1 l.18[P.1 l.18] ``h'' added\\ 
P.6 l-14 to l-1[P.6 l-12\&-13] two sentences removed, became new paragraph 
starting at P.6 l.-14\\
P.7 l.1[P.6 l.-5] ``$k \neq 0$'' replaced with ``$k\in \EuScript{S}$''\\
P.7 l.4[P.6 l.-2] ``$k \neq 0$'' replaced with ``$k\in \EuScript{S}$''\\
P.11 l.-3[P.11 l.-8] ``decay'' replaced with ``growth''\\
P.12 l.15[P.12 l.15] ``at most'' added\\ 
P.15 l.24[P.15 l.24] reference [Oel-Zaf] updated\\
}

\begin{abstract}
We consider a family $\big\{ E_m(D,M) \big\}$ of holomorphic bundles
constructed as follows:
to any given $M\in GL_n(\Z)$, we associate a ``multiplicative automorphism''
$\varphi$ of $(\C^*)^n$. Now let $D\subseteq (\C^*)^n$ be a $\varphi$-invariant
Stein Reinhardt domain. Then $E_m(D,M)$ is defined as the flat bundle
over the annulus of modulus $m>0$, with fiber $D$, and monodromy $\varphi$.

\smallskip
We show that the function theory on $E_m(D,M)$ depends nontrivially
on the parameters $m, M$ and $D$. Our main result is that
$$E_m(D,M) \text{\ is Stein
if and only if\ } m \log \rho (M) \leq 2 \pi^2,$$
where $\rho(M)$ denotes
the max of the spectral radii of $M$ and $M^{-1}$.

\smallskip
As corollaries, we:\\
-- obtain a classification result for Reinhardt domains in all dimensions;\\
-- establish a similarity between two known counterexamples to a question
of J.-P. Serre;\\
-- suggest a potential reformulation of a disproved conjecture of Siu Y.-T.
\end{abstract}
\maketitle \thispagestyle{empty}

\noindent
Let $D$ be a Stein manifold. We say that $D$ belongs to $\S$ when:
for any Stein manifold $B$ and any locally trivial bundle $E \rightarrow B$,
the manifold $E$ is also Stein. A famous question of Serre can be formulated as:
``Are all manifolds in \S?"

Skoda answered it in the negative, by proving that $\C^2 \not\in \S$ (\cf \cite{Sko}).
Mok showed that any open Riemann surface belongs to $\S$ (\cf \cite{Mok}).

Many bounded domains belong to $\S$: for any given bounded domain $D$,
Diederich-Forn\ae ss-Stehl\'e showed that if $\partial D$
is smooth, then $D\in \S$ (\cf \cite{Die-For} and \cite{Steh}).
Siu showed that if $b_1(D)=0$ then $D\in\S$ (\cf \cite{Siu}).

However, C\oe ur\'e and L\oe b ({\em cf.} \cite{Coe-Loe}) found a bounded
domain $D_{CL}$ not in $\S$. It happens that $D_{CL}$ has the Reinhardt symmetry.

It is an open problem to characterize bounded domains of $\C^d$ not in $\S$ (\cf \cite{CZ}).
This classification problem is solved for all bounded {\em Reinhardt}
domains with $d=2$ in \cite{Pfl-Zwo} and with $d=2$ or $3$ in \cite{OZ}.

\smallskip
Here we study holomorp\modifinred{h}ic functions on a family of bundles $\{E_m(D,M)\}$ over annuli,
depending on a non necessarily bounded Reinhardt domain $D\subset \Csd$, a matrix $M$ and a number $m$.
This family contains the non-Stein bundle of C\oe ur\'e and L\oe b ({\em cf.} \cite{Coe-Loe}).
Moreover, roughly speaking, all known examples
of bounded domains not in $\S$ appear as fibers in these bundles.
Our main motivations are the problems of characterizing bounded (Reinhardt) domains not in $\S$,
and also, when $D\not\in \S$ is given, of characterizing the Steinness of a bundle with fiber $D$.

\section{Setting and results}\noindent
We work in the category of complex manifolds and holomorphic maps.

\medskip
Let $d\geq 2$ and $M\in \gldz$. Then $M$ gives an automorphism of $\Csd$ defined by
\mbox{$\z=(z_1,\dots, z_d)\mapsto (w_1,\dots,w_d)$} with $w_i= z_1^{M_{i1}}\dots z_d^{M_{id}}$.
This automorphism induces a \Z-action. We denote by $j.\z$ the action of $j\in \Z$ on $\z\in\Csd$.
Fix $D\subset \Csd$ any Stein \Z-invariant Reinhardt domain. \cmt{Paragraph slightly re-arranged}

\smallskip
Let $m>0$. We realize the annulus of modulus $m$ by
$$A_m=\setofst{w\in\C}{1<|w|<e^m},$$
and its universal cover by the horizontal strip
$$S_m=\left\{w\in \C\ \Big\vert \quad |\mathit{Im}\ w|< \frac{m}{4 \pi} \right\}$$
endowed with the \Z-action generated by $w\mapsto w+ 1$.
Then $S_m / \Z$ is isomorphic to $A_m$.

\smallskip

There is a diagonal \Z-action on $S_m \times D$ given by $j.(w,\z)=(w+ j,j.\z)$,
which is free and properly discontinuous. We consider the quotient manifold
$$E_m=E_m(D,M)=\frac{S_m \times D}{\Z}.$$

The projection map $S_m \times D \rightarrow D$ induces a locally trivial fibration
$$p: E_m \rightarrow A_m$$ with fiber $D$. \cmt{Notation explained}
We extend the notation to cover the cases
of the two annuli of infinite modulus: $A_{\infty}=\setofst{w\in\C}{1<|w|<\infty}$ ($=$the pointed disk)
and $A_{2\infty}=\setofst{w\in\C}{0<|w|<\infty}$ ($=\Cs$), with corresponding $S_{\infty}$ (=\H),
$S_{2\infty}$ ($=\C$) and bundles $E_{\infty}$, $E_{2\infty}$.

\smallskip
We denote by $\rho(M)$ the max of the spectral radii of $M$ and $M\inv$.
\medskip
\begin{MainTheorem}
$E_m(D,M)$ is Stein if and only if $\ m\, \log \rho(M)\leq 2 \pi^2$.
\end{MainTheorem}

\medskip

\noindent
Notice that we do not assume boundedness of $D$.
By the methods of \cite{OZ}, one can prove the Steinness when $\rho(M)=1$.
The results in that paper also imply the ``only if" part for $d=2$, and
implicitly for $d=3$.

The Steinness when $\rho(M)>1$ and $m$ is small enough is new even in small dimensions.
Sections \ref{only_if} and \ref{if} are devoted
respectively to the proofs of the ``only if" and ``if" parts of the Main Theorem.
They both rely upon the existence or non-existence of fast decaying functions on $S_m$,
and a Laurent series method already used in \cite{Zaf}.

\subsection{Case of a bounded $D$.}
\label{case_bounded_D}


The ``only if" part of the Main Theorem leads 
to a classification result that we now explain.

Let $D\subset\Csd$ be a bounded Reinhardt domain.
Theorem 3 in \cite{OZ} yields: if no $M$ with $\rho(M)>1$
multiplicatively acts on $D$ (see below), then $D\in\S$.

\medskip
Let $m>0$.
Let $E'\rightarrow A_m$ be {\em any} fiber bundle over the annulus $A_m$
with fiber $D$.
It follows from \cite{Roy} that $E'$ is flat, so the monodromy
$h\in\mathit{Aut}(D)$ characterizes $E'\rightarrow A_m$.
By \cite{Shim} we know that $h$ is of the form
$$h(z_1,\dots,z_d)=
(c_1\; z_1^{M_{11}}\!\!\!\dots z_d^{M_{1d}},\dots\dots, c_d\; z_1^{M_{d1}}\!\!\!\dots z_d^{M_{dd}}),$$
with $M\in GL_d(\Z)$ and $c_1,\dots,c_d \in \Cs$.
When a given $D$ admits such an automorphism, we say that
$M$ {\em multiplicatively acts} on $D$.

The bundle $E'$ is very similar to $E_m(D,M)$.
In fact it follows from the Main Theorem that
$E'$ is Stein if and only if $m\log \rho(M)\leq 2\pi^2$: \cmt{One paragraph moved from previous Sect.4 to here}
If $\mathit{Spec}\ M=\{1\}$ (so $\rho(M)=1$), then by Theorem 3 in \cite{OZ},
$E'$ is Stein. Up to applying (several times) Lemma 1.7 \cite{OZ}, we can assume that
$1\not\in \mathit{Spec}\ M$. Then, as shown in 3.3 \cite{OZ},
there is an automorphism of $\Csd$ that
conjugates $h$ to
$g:(z_1,\dots,z_d)\mapsto
(z_1^{M_{11}}\!\!\!\dots z_d^{M_{1d}},\dots\dots, z_1^{M_{d1}}\!\!\!\dots z_d^{M_{dd}})$.
Therefore $E'$ is isomorphic to $E_m(D,M)$ and we apply the Main Theorem.

The above implies in particular that as soon as $M$ with $\rho(M)>1$
multiplicatively acts on $D$, there exists a non-Stein bundle with fiber $D$ over a Stein base
(a thick enough annulus). So $D\not\in\S$.

Summing up, the Main Theorem together with Theorem 3 in \cite{OZ} imply the
\begin{thm}
\label{which_D_s_are_in_S_dim_d}
Let $D\subset (\Cstar)^d$ be a bounded Stein Reinhardt domain, with $d\geq 2$.
Then $D\in \S$ 
if and only if no $M$ with $\rho(M)>1$
multiplicatively acts on $D$.
\end{thm}\cmt{One paragraph removed}
\noindent
From the known cases of $d=2$ or $3$,
we expect that the assumption that $D$ do not intersect
any coordinate hyperplane is only technical, and the result should hold
unchanged for any bounded Reinhardt domain $D \subset \C^d$.

We use here a different method than in \cite{OZ} to obtain
the higher-dimensional case, but
the trade-off is that we lose some information.
We don't know which matrices can appear for a bounded $D$ not in $\S$, thus
we don't know much about the geometry of such domains.
However, the situation {\em does} become more complicated in dimension
at least $4$.
For $d=2,3$, we showed in \cite{OZ} that a bounded $D$ not in $\S$
can only be multiplicatively acted on by a real-diagonalizable matrix. This
leads to a relatively simple geometric description of these domains.
These facts do not extend to $d\geq 4$:

\inred{
\begin{prop}
There exists a Stein bounded Reinhardt domain $D$ in $(\Cs)^4$
that is multiplicatively acted on by a matrix $M\in GL_4^{\rho > 1}(\Z)$ with non-real spectrum.
In particular, by Theorem \ref{which_D_s_are_in_S_dim_d}, $D\not\in\S$.
\end{prop}
      }
\begin{proof}
Define
$$N:=
\begin{bmatrix}
0& -2& -7& 9\\
0&-10&-20&29\\
0&-13&-31&43\\
-1&-11&-36&47
\end{bmatrix}
\in GL_4(\Z).$$
Then $Spec\ N=\{ \alpha_1= 6.23...,\ \alpha_2= 0.27...,\ \omega= -0.25...\!\!+ i 0.73...,\ \bar{\omega} \}$.
In particular, $N$ is not similar to a block-diagonal matrix in $GL_2(\Z)\oplus GL_2(\Z)$.

Direct computation shows that $N$ admits real eigenvectors $v_1$ and $v_2$ associated to $\alpha_1$ and $\alpha_2$,
and such that
\def\rnegf{(\R^{\scriptscriptstyle{<0}})^{4}}
$v_i \in \rnegf, i=1\dots 2$. Let $\{w', w''\}$ be a basis of the $N$-invariant real two-dimensional subspace
of $\R^4$ corresponding to the complex eigenvalue $\omega$.

Now take $u:=[-14,-43,-62,-63]^t$. 
Direct computation shows that $u$ decomposes as
$u=a_1 v_1 +a_2 v_2 + a' w' + a'' w''$ with $a_1>0$ and $a_2>0$.

Denote the iterates $u_j:= N^j u= (\alpha_1)^j a_1 v_1 +(\alpha_2)^j a_2 v_2 + w_j$.
Then $\norme{w_j}$ is of order $|\omega|^j$, but $|\alpha_2| < |\omega| < |\alpha_1|$
and $|\alpha_2| < 1 < |\alpha_1|$, so there exists $J>0$ such that: for all $j\in\Z$,
if $|j|\geq J-4$ then $u_j \in \rnegf$.

Define $A:=\{u_{(2j+1)J+k}\ |\ j\in\Z, k=0\dots 4 \}$, and $\Omega$ as the interior
of the convex hull of $A$ in $\R^4$. Then
\begin{itemize}
\item $\Omega$ is $N^{2J}$-invariant, because $A$ is;
\item $\Omega$ is contained in $\rnegf$, because $|(2j+1)J+k|\geq J-4$;
\item $\Omega$ is not empty:
By a rank computation (over \Z), one checks that $B:=\{u_k\ |\ k=0\dots 4\}$ is not contained
in any affine hyperplane of $\R^4$. Then, the same is true of $A$ because
$A$ contains $\{u_{(2j+1)J+k}\ |\ j=0, k=0\dots 4 \}=N^J B$.
\end{itemize}
Then $D:=\{ \z \in  (\Cs)^4\ |\ \ (\log |z_1|,\dots, \log |z_4|) \in \Omega \}$ is a non-empty, bounded,
Stein Reinhardt domain multiplicatively acted on by $M=N^{2J}$, and $\rho(M)= (\alpha_1)^{2J}>1$.
\end{proof}


\subsection{Case of $D=\Csd$.} Take $M\in GL_2(\Z)$ with $\rho(M)>1$.
Take $E_{2\infty}=E_{2\infty}\big( (\Cs)^2, M \big)$, with fibration map
$p:E_{2\infty}\rightarrow \Cs$. Then ({\em cf.} \cite{Coe-Loe})
not only is $E_{2\infty}$ not Stein, but moreover
$p^*:\O(\Cs) \rightarrow \O(E_{2\infty})$ is an isomorphism, i.e.,
all functions on $E_{2\infty}$
come from the base.
In any dimension and over any annulus we prove the
\begin{thm}
\label{fiber_Csd}
Take $d\geq 2$, $M\in\gldzrbo$ and $m$ any positive number.
Let $E_m=E_m \big( \Csd,M \big)$.
If the characteristic polynomial
of $M$ is irreducible over $\Z$ and $m \log \rho(M)>2\pi^2$,
then the fibration map $p:E_m\rightarrow A_m$ induces an isomorphism between
functions on $A_m$ and $E_m$. (Proof is in Sect. \ref{proof_of_fiber_Csd})
\end{thm}
\noindent
So for any bundle $E_m$ considered in this theorem, there is a
critical value for $m$ below which $E_m$ is Stein, and above which
the only functions come from the base.
As opposed to this situation,
recall that when $D$ is a bounded domain and $E_m$ is a bundle
with fiber $D$ over any Stein base, it is known
that functions on $E_m$ separate points ({\em cf.} \cite{Siu}).

\subsection{Remark: Connection with the Schinzel-Zassenhaus problem.}
\label{SectionSZ}
\cmt{Sect. \ref{SectionSZ} rewritten}
From the ``only if" part of the Main Theorem and Theorem \ref{fiber_Csd}, we obtain the
\begin{crl}
\label{mu}
There exists a function $\mu:\N\rightarrow \R^{>0}$ such that for all
$d>0$ and $M\in\gldzrbo$:
\begin{enumerate}
\item[(a)] for any Reinhardt domain $D$ invariant by the multiplicative action of $M$, if $m>\mu(d)$ then $E_m(D,M)$ is not Stein;
\item[(b)] if $m>\mu(d)$ and the characteristic polynomial of $M$ is irreducible, then any function on $E_m\big(\Csd,M\big)$ is constant on fibers.
\end{enumerate}
\end{crl}
\begin{proof}
For any polynomial $P$, denote by $\rho(P)$ the maximal modulus of its roots. Fix an integer $d\geq 1$.
From the relations between the roots and coefficients of a polynomial,
it follows that for any $a\in\R$, the set
$$\setofst{P\in\Z[X]}{P \text{ monic},\ \deg P=d,\ \rho(P)\leq a}$$
is finite. The corresponding set of roots is finite, so there exists
$\mu'(d)>0$ such that if $P$ is any monic, integral polynomial of degree $d$,
$\rho(P)>1$ implies $\rho(P)\geq 1+\mu'(d)$. We take $\mu'(d)$ minimal with that property.

Now let $M\in\gldzrbo$ with characteristic polynomial $P_M$.
Up to inverting $M$, we can assume that $\rho(P_M)=\rho(M)$.
By assumption, $\rho(M)>1$, so $\rho(M)\geq 1+\mu'(d)$. Thus if we
define $\mu(d):=\frac{2\pi^2}{\log\big( 1+\mu'(d) \big)}$, the results follow
from the Main Theorem and Theorem \ref{fiber_Csd}.
\end{proof}

It is of independent interest to find estimates of $\mu'(d)$.
Note that $\mu'(d) \leq \rho(X^d-2) -1 \sim \frac{\log 2}{d}$. 
In 1965, Schinzel and Zassenhaus asked in \cite{SZ} whether there exists $\gamma>0$ independent
of $d$ such that $\mu'(d) \geq \frac{\gamma}{d}$.

Only partial results are known (see the survey \cite{Smy}).
In 1971, Smyth answered that question positively for the class of all
non-reciprocal polynomials.
One of the best general results, proved by Voutier in \cite{Vou}, is
$$\mu'(d) \geq \frac{1}{4d}\left(\frac{\log\log d}{\log d}\right)^3.$$
Each of these estimates yields an estimate on the function
$\mu(d)=\frac{2\pi^2}{\log\big( 1+\mu'(d) \big)}$.

On the other hand, if $P$ is a monic, integral, degree $d$ reciprocal polynomial,
its companion matrix belongs to $\gldz$.
Therefore, by the ``if" part of the Main Theorem and Smyth's result, the Schinzel-Zassenhaus problem is
equivalent to the existence of $\gamma$ such that:
$$\text{If } m>\frac{2\pi^2}{\log(1+\frac{\gamma}{d})}
\text{ then } \forall M\in\gldzrbo,\ E_m \text{ is not Stein.}$$
Note that asymptotically with respect to $d$, the condition 
is simply
$m>\gamma' d$. 

\subsection{Remark: Stein bundles over a non-Stein base.}
Let $E \rightarrow B$ be a fiber bundle with fiber $D$. As $D$ appears as
a closed submanifold of $E$, it is necessary that $D$ be Stein for $E$ to be Stein.
On the other hand, it can happen that $B$ is not Stein but $E$ is
(e.g., $SL_2(\C)\rightarrow \C^2-\{0\}$ with fiber $\C$).
We refer to the
articles by M. Abe for related results, and simply notice the following facts.
Take $d=2$ and $M\in\gldzrbo$. Then for certain choices of a bounded $D$ (e.g., $D=D_{CL}$),
the manifold $E_m$ admits {\em another} fibration $E_m\rightarrow \cusp$ with fiber
$S_m$ and base a non-Stein manifold $\cusp$ studied by Hirzebruch
({\em cf.} \cite{Zaf}). Note that the pair $(S_m, w\mapsto w+1)$
is isomorphic to $(\Delta, \psi)$,
where $\Delta$ is the unit disk in $\C$ and $\psi$ some automorphism.
Varying $m$ corresponds to varying $\psi$.
Then one can construct a continuous family $\{B_t\}_{t\geq 0}$ of
disk bundles over $\cusp$ such that:
$B_0$ is the trivial bundle (corresponding to $\psi$ being the identity);
by the Main Theorem,
there exists $s>0$ such that $B_t$ is Stein if and only if $t\geq s$.

\section{Main Proofs.}\label{proofs}

\subsection{Proof of the ``only if" part of the Main Theorem.}
\label{only_if}

\noindent
We want to show: If $E_m$ is Stein then $\ m\, \log \rho(M)\leq 2 \pi^2$.

For the sake of clarity, we denote vector quantities by bold characters.
We denote by $$\applic{q}{S_m \times D}{\ \EMm}{(w,\z)}{\ [w,\z]}$$
the quotient map that was used to define $\EMm$. It induces an
isomorphism between the spaces of scalar-valued functions $\O\big(\EMm\big)$
and $\O^\Z\big(S_m \times \Csd\big)$.

For $f\in\O\big(\EMm\big)$, we will use implicitly that isomorphism and
write indifferently $f[w,\z]$ or $f(w,\z)$.

We denote an element $\k \in \Z^d$ by a row vector $\k=(k_1,\dots,k_d)$.
For $\z \in \Csd$, we denote by $\z^{\sk}$ the product $z_1^{k_1}\dots z_d^{k_d}\ \in \C$.

Let $\Delta \subset S_m$ be a small disk \inred{centered at $w_0$}.
Then $\Delta\times D$
is a Reinhardt domain, so we can expand $f$ there into a Laurent series:
\begin{equation}
\label{local_Laurent}
\forall (w,\z)\in \Delta \times D, \quad
f(w,\z)= \sum_{i\in\N, \sk\in\Z^d} a_{i\sk} (w-w_0)^i \z^\sk.
\end{equation}
Thanks to absolute convergence of such a series, we can write
\begin{equation}
\label{grouped_coeffs}
f(w,\z)
= \sum_{\sk} \Big( \sum_i a_{i\sk} \inred{(w-w_0)}^i\Big) \z^\sk
= \sum_{\sk} g_\sk(w) \z^\sk.
\end{equation}
Notice that by the uniqueness of Laurent expansions, for all $\k$, the function $g_\sk$
so defined is actually
a well-defined element of $\O(S_m)$, \inred{and (\ref{grouped_coeffs}) is valid on $S_m \times D$}. Moreover,
as (\ref{local_Laurent}) is locally absolutely uniformly converging with
respect to $(w,\z)$,
so is the rightmost series in (\ref{grouped_coeffs}).
This expansion of $f$ is called a Hartogs-Laurent
series.

For $j\in \Z$ and $\k \in \Z^d$, we denote $j.\k=\k M^j$. 
\modifinred{This defines a \Z-action on $\Z^d$ and moreover,} 
for all $j$, $\z$ and $\k$, 
\begin{equation}
\label{jzk_equals_zjk}
(j.\z)^{\sk}=\z^{j.\sk}.
\end{equation}
As $f$ is $\Z$-invariant, we know that for all $j,w$ and $\z$, 
$f(w,\z)=f(w+j,j.\z)$,
so by (\ref{grouped_coeffs})
and (\ref{jzk_equals_zjk}),
$$\sum_{\sk} g_\sk(w) \z^\sk = \sum_{\sk} g_\sk(w+j) \z^{j.\sk}.$$
Uniqueness of the Hartogs-Laurent expansion implies that for all $j,\k$ and\nolinebreak\ $w$,
$g_{j.\sk}(w)=g_\sk(w+j)$. 
\def\SctnZero{\EuScript{S}_0}

\modifinred{
Let $\SctnZero$ be the set of elements in $\Z^d$ with a finite \Z-orbit. 
If $\SctnZero = \Z^d$, then $\mathit{Spec}\ M \subset S^1$, so 
$\rho(M)=1$; thus the ``only if'' part of the theorem becomes trivially true. Hence 
we can assume that $\Z^d - \SctnZero$ is non-empty. The \Z-action restricted to 
$\Z^d - \SctnZero$ is free. Choosing an arbitrary section $\EuScript{S}$ of this free action, 
we can write 
$$f(w,\z)= \sum_{\sk\in\SctnZero}  g_\sk(w) \z^{\sk} 
         + \sum_{\sk\in\EuScript{S}} \sum_{j\in\Z}  g_\sk(w+j) \z^{j.\sk}.$$}
Now for any fixed $\k\in \EuScript{S}$, the series
$$
\sum_{j\in\Z}  g_\sk(w+j) \z^{j.\sk}
$$
is a subseries of (\ref{grouped_coeffs}) because $g_\sk(w+j)=g_{j.\sk}(w)$,
hence it is locally absolutely uniformly converging.
In particular, for all $\z \in D$ and 
\modifinred{$\k\in \EuScript{S}$}
\begin{equation}
\label{fast_decay}
g_\sk(w+j) \z^{j.\sk} \xrightarrow{j\rightarrow \pm\infty} 0
\end{equation}
locally uniformly with respect to $w$.

\medskip
{\bf Vocabulary.} A subset $J\subset \Z$ is said
to have {\em bounded gaps} not bigger than $l$ when:
for all $j\in J$, there exists $j'\in J-\{j\}$ such that $|j-j'| \leq l$.



\medskip
The proof of the following proposition, which is a Phragm\'en-Lindel\"of-type result 
(\cf \cite{Tit} 5.65), 
is directly adapted from a proof 
that was kindly communicated to me by A. Baernstein \inred{and L.
Kovalev}.
\begin{prop} 
\label{P-L}
Let $\delta, \mu,h \in \R$ such that $\delta >0$, $\mu >0$ and $h>\frac{\pi}{\mu}$.
Let $S=\setofst{w\in\C}{|\mathit{Im}\ w|<\frac{h}{2}}$ and $g\in \O(S)$.
Assume that $w$-locally uniformly,
$$g(w+j)\ e^{\delta e^{j \mu_+}} \xrightarrow[j\in J_+]{j\rightarrow +\infty} 0
\quad \text{ and } \quad
g(w+j)\ e^{\delta e^{|j|\mu_-}} \xrightarrow[j\in J_-]{j\rightarrow -\infty} 0,$$
where $J_+$,$J_-$ are infinite subsets of\, $\Z$ with bounded gaps,
and 
$\max(\mu_+, \mu_-)= \mu$.
Then $g$ vanishes identically.
\end{prop}
\begin{proof}
Choose an $h'$ such that $\frac{\pi}{\mu}<h'<h$.
For $\alpha \in (0,+\infty]$, define
$$R_\alpha = \setofstquad{w\in\C}{ |\mathit{Re}\ w|\leq \alpha\ ,\ |\mathit{Im}\ w|\leq h'},$$
$R^+_\alpha=R_\alpha \cap \{ \mathit{Re}\ w \geq 0 \}$ and
$R^-_\alpha=R_\alpha \cap \{ \mathit{Re}\ w \leq 0 \}$.

Let $l$ be a majorant of the gaps of $J_+$ and $J_-$.
By assumption there exists $j_0 > l$ such that:
$$\forall j\in J_+\cap (j_0, +\infty), \ \forall w\in R_l +j,
\qquad \log |g(w)| \leq -\delta e^{j \mu_+}, \text{and}$$
$$\forall j\in J_-\cap (-\infty, -j_0), \ \forall w\in R_l +j,
\qquad \log |g(w)| \leq -\delta e^{|j|\mu_- }.$$
On the other hand,
$$\forall j\in J_+\cap (j_0, +\infty),\ \forall w\in R_l^- +j,
\qquad |\mathit{Re}\ w| <j,$$
$$\forall j\in J_-\cap (-\infty, -j_0),\ \forall w\in R_l^+ +j,
\qquad |\mathit{Re}\ w| <|j|.$$
As the rectangles $R_l^{\pm}$ are wider than the gaps of $J^{\pm}$, we get
$\forall w\in R_\infty$,
$$\text{if } \mathit{Re}\ w > j_0  \text{ then } \log |g(w)| \leq -\delta e^{\mu_+ \mathit{Re}\ w}, \text{ and}$$
$$\text{if } \mathit{Re}\ w < -j_0 \text{ then } \log |g(w)| \leq -\delta e^{\mu_- |\mathit{Re}\ w|}.$$
Define $h\in\O(S)$ by $h(w)=g(w) g(-w)$.
Then $\log |h(w)|= \log |g(w)| + \log |g(-w)|$, so by the previous inequalities,
$$\forall w\in R_\infty -R_{j_0},\ \log |h(w)| \leq -\delta e^{\mu |\mathit{Re}\ w|}.$$

In particular, $\sup_{\partial R_\infty}\ \log |h| = C_1 < +\infty$. Take any $\alpha> j_0$.
Let $\omega$ be the harmonic function in the interior of the rectangle $R_\alpha$ with
boundary values $0$ on horizontal sides and $1$ on vertical sides.
One can find an explicit expression of $\omega$ as a series of functions, from
which one can get a $C_2>0$ independent from $\alpha$ such that
$$\omega(0)\geq C_2 e^{-\frac{\pi\alpha}{h'}}.$$
On the boundary of $R_\alpha$, the subharmonic function
$w\mapsto \log |h(w)|$ is majorized by the harmonic function
$w\mapsto C_1 - \delta e^{\mu\alpha} \omega(w)$.
Therefore 
$$\log |h(0)| \leq
C_1 - \delta e^{\mu\alpha} C_2 e^{-\frac{\pi\alpha}{h'}} =
C_1 - \delta C_2 e^{(\mu-\frac{\pi}{h'})\alpha}.$$
From $h'> \frac{\pi}{\mu}$ we get $\mu - \frac{\pi}{h'} >0$.
As the inequality is true for all $\alpha > j_0$, we conclude that $h(0)=0$.

Now, for any $\varepsilon$ so small that $R_\alpha  + i \varepsilon \subset S$,
we can repeat the above argument and obtain that $h(i \varepsilon)=0$.
By the principle of isolated zeros, $h$ vanishes identically on $S$.
Again by the principle of isolated zeros, it follows that $g$ also vanishes identically.
\end{proof}

We denote $$\mathbf{log}\ |\z|=\vect{\log |z_1| \\ \vdots \\ \log |z_d|}\in \R^d.$$
Then for all $\k$ and $\z$, a direct computation gives
\begin{equation}
\label{ln_zk}
\log |\z^{\sk}|=\innerp {\k}{\mathbf{log}\ |\z|}= k_1 \log |z_1|+ \dots + k_d \log |z_d|\ \in \R,
\end{equation}
where
$\innerp{{\scriptstyle \bullet}\ }{{\scriptstyle \bullet}}$ denotes the usual inner product in $\R^d$.

For $p\geq 1$, denote $\T=(S^1)^p \subset \C^p$.
For $\theta=(\theta_1,\dots, \theta_p) \in \T$ and $j\in \Z$,
denote $\theta^j=(\theta_1^j,\dots, \theta_p^j)$. We will use the elementary
\begin{lem} 
\label{bounded_gaps}
Let $\varepsilon>0$. For any $\theta \in \T$, 
the set $$J=\{j\in \N\ |\quad \norme{(1,\dots,1) - \theta^j}_\infty <\varepsilon \}$$
is infinite and has bounded gaps.
\end{lem}
\begin{proof}
Consider the semigroup $X=\setofst{\theta^j}{j\in\N} \subset \T$.
Its closure $\bar{X}$ is a closed semigroup in a compact group, hence is a group.

If $\bar{X}$ is discrete then all $\theta_i$'s are roots of unity and the conclusion
follows easily, so we can assume that $\bar{X}$ is a subtorus
of $\T$ of dimension $n>0$. For simplicity we assume that $n=2$
(the case $n=1$ is easier, and $n\geq 3$ is similar).

Denote $V=\text{Lie}(\bar{X})$, which is a $2$-dimensional subspace of $\text{Lie}(\T)\approx \R^p$.
Denote $U(\varepsilon)$ the $\varepsilon$-neighborhood of $(1,\dots,1)$ in $\T$.
There exist a small neighborhood $A$ of $0$ in $V$ and a small neighborhood $B$ of
$(1,\dots,1)$ in $\bar{X}$ such that $B\subset U(\varepsilon)$ and
$\mathit{exp}: A\rightarrow B$ is a local group isomorphism.
Take linear coordinates $v_1, v_2$ in $V$. We can assume that $A$
is a ball for the associated sup-norm. Consider the open quadrants
$Q_1,\dots, Q_4 \subset V$ determined by these coordinates.

As $X$ is dense in $B$, we can find $j_1,\dots, j_4 \in \N$ such that
each quadrant contains exactly one element of
$\{ \alpha_1=\mathit{exp}\inv(\theta^{j_1}), \dots , \alpha_4=\mathit{exp}\inv(\theta^{j_4}) \}$.
Now define a sequence $(a_p)_{p\in \N}$ in $A$ by:
$a_0= \alpha_1$, and $a_{p+1}= a_p+ \alpha_i$, where $\alpha_i$ is contained in
the quadrant opposite to any closed quadrant containing $a_p$.
Now the images by $\mathit{exp}$ of this sequence give a subsequence
of $(\theta^j)_{j\in \N}$ contained in $B\subset U(\varepsilon)$
showing that $J$ is infinite with gaps not bigger than $\max \{j_1, \dots,j_4 \}$.
\end{proof}

Denote $\mu = \log \rho(M)$. If $\mu =0 $ there is nothing to
prove, so we assume $\mu>0$. We choose a numbering of the
eigenvalues of $M$ such that
$$|\lambda_1|=\dots= |\lambda_s| >
|\lambda_{s+1}| \geq \dots \geq |\lambda_{t-1}| >
|\lambda_t|=\dots = |\lambda_d|.$$
By our assumption on $\mu$,
$|\lambda_{1}|
> 1 > |\lambda_d| > 0$. Denote $\mu_+ = \log |\lambda_{1}|$ and
$\mu_- = - \log |\lambda_{d}|$. Then by definition of $\rho(M)$,
$\mu= \max(\mu_+ , \mu_-)$.

\begin{lem} 
\label{zd_and_kd}
There exist \emph{\zd}$\in D$ and \emph{\kd}$\in \Z^d$ such that for some $\delta>0$,
\emph{$|\zd^{j.\skd}| > e^{\delta e^{j \mu_+}}$} for all $j\in J_+$, and
\emph{$|\zd^{j.\skd}| > e^{\delta e^{|j| \mu_-}}$} for all $j\in J_-$,
where $J_+\subset \Z^{\geq 0}$ and $J_-\subset \Z^{\leq 0}$ are
infinite with bounded gaps.
\end{lem}
\begin{proof}
Let $\{ \bo, \dots, \bd \}\subset \C^d$ be a Jordan basis for the linear
action of $M$ on $\C^d$.
Up to exchanging some eigenvalues of equal modulus, we can assume that this
numbering matches the above numbering of the eigenvalues.

Pick $\zd \in D$, and denote
$\mathbf{log}\ |\zd| =\ub
= \sum_{i} u_i \bi
\in \R^d \subset \C^d$.
As $D$ is open we can assume that $u_1 \neq 0$ and $u_d \neq 0$.

When $j\rightarrow +\infty$, the dominant terms
of $M^j\ub$ come from the diagonals of the first blocks in the Jordan form of $M$.
Namely we can write
$$
M^j\ub = |\lambda_1|^j \Big( \underbrace{\sum_{i=1}^{s} \theta_i^j u_i \bi }_{\ub'_{+,j}} + \ub''_{+,j} \Big)
$$
with $\ub''_{+,j} \xrightarrow{j\rightarrow +\infty} 0$ and
$\theta_i = \frac{\lambda_i}{|\lambda_1|} \in S^1, \ i=1\dots s$.

And similarly
$$
M^j\ub = |\lambda_d|^j \Big( \underbrace{\sum_{i=t}^{d} \theta_i^j u_i \bi }_{\ub'_{-,j}} + \ub''_{-,j} \Big)
$$
with $\ub''_{-,j} \xrightarrow{j\rightarrow -\infty} 0$ and
$\theta_i = \frac{\lambda_i}{|\lambda_d|} \in S^1, \ i=t\dots d$.

Denote $\ub'_+= {\displaystyle \sum_{i=1}^{s} u_i \bi }$ and
$\ub'_-= {\displaystyle \sum_{i=t}^{d} u_i \bi }$.
Notice that $\ub'_+$ and $\ub'_-$ are linearly independent, and real
because $M$ and $\ub$ are real.
In particular, we can pick $\kd \in \Z^d$ such that
$$
\innerp{\ub'_+}{\kd}\ > 3\delta  \text{ and } \innerp{\ub'_-}{\kd}\ > 3\delta
$$
for some $\delta>0$.
By applying Lemma \ref{bounded_gaps} to $\theta=(\theta_1, \dots ,\theta_s)$
for a small enough $\varepsilon$, we get an infinite
$J_+ \subset \Z^{{\scriptscriptstyle\geq 0}}$ with
bounded gaps, such that $$\forall j\in J_+,\ \innerp{\ub'_{+,j}}{\kd}\ > 2\delta .$$
By applying Lemma \ref{bounded_gaps} to $\theta=(\theta_t\inv, \dots ,\theta_d\inv)$
for a small enough $\varepsilon$, we get an infinite
$J_- \subset \Z^{{\scriptscriptstyle \leq 0}}$ with
bounded gaps, such that $$\forall j\in J_-,\ \innerp{\ub'_{-,j}}{\kd}\ > 2\delta .$$

\smallskip
Now, by (\ref{ln_zk}) and (\ref{jzk_equals_zjk}), for all $j$, $\k$ and $\z$,
$$
\log |\z^{j.\sk}|=\innerp {\mathbf{log}\ |\z|}{\k M^j} = \innerp{ M^j(\mathbf{log}\ |\z|)} {\k}.
$$
Thus for all
$j\in J_+,\ \log |\zd^{j.\skd}| = |\lambda_1|^j(\innerp{\ub'_{+,j}}{\kd} + \innerp{\ub''_{+,j}}{\kd})$,
and as $\ub''_{+,j} \xrightarrow{j\rightarrow +\infty} 0$, by discarding a finite number of (not big enough) elements from $J_+$, we get
$$\forall
j\in J_+,\ \log |\zd^{j.\skd}| > \delta |\lambda_1|^j = \delta e^{j\mu_+}.$$
Similarly, by discarding a finite number of elements from $J_-$, we get
$$\forall
j\in J_-,\ \log |\zd^{j.\skd}| > \delta |\lambda_d|^{-|j|} = \delta e^{|j|\mu_-}.$$
\end{proof}

We now show by contradiction that $m\, \log \rho(M)\leq 2 \pi^2$:
Assume that $m > \frac{2 \pi^2}{\log \rho(M)}$.
By (\ref{fast_decay}), Proposition \ref{P-L} and Lemma \ref{zd_and_kd}, it follows that $g_{\skd}$
vanishes identically. But $\EMm$ is Stein, so any holomorphic function defined
on some fiber (which is a closed submanifold) extends to a function on $\EMm$.
Equivalently, any function on $\{w\}\times D$ extends to a $\Z$-invariant
function on $S_m \times D$. In particular, there exists such a function $f$
satisfying $f(0,\z)= \z^{\skd}$. This function's Hartogs-Laurent expansion
(\ref{grouped_coeffs}) has a non-zero coefficient corresponding to $\z^{\skd}$.
So $g_{\skd}$ can
not vanish identically.

\subsection{Proof of the ``if" part of the Main Theorem.}
\label{if}

\noindent
We want to prove: If $\ m\, \log \rho(M)\leq 2 \pi^2$ then $E_m$ is Stein.

By construction $D\hookrightarrow\Csd$ is $\Z$-equivariant, so
$E_m(D,M)$ can be seen as subbundle of $E_m\big(\Csd,M\big)$.
In particular, $E_m(D,M)$ is a locally Stein open subset of $E_m\big(\Csd,M\big)$.
By the Docquier-Grauert
theorem, it is therefore enough to prove the statement under the assumption
that the fiber $D$ is $\Csd$. We will simply write $E_m$ instead
of $E_m\big(\Csd,M\big)$.

{\em Outline of the proof:} we roughly follow the line of argument of \cite{Siu}.
The key result \inred{here} is Lemma \ref{any_monomial}, which
plays the role of Siu's ``Main Lemma". For
this result the changes are of course essential
because neither of the hypotheses made in \cite{Siu}
(vanishing of $b_1(D)$ and boundedness of $D$) holds here.

First we will assume that $m\log \rho(M)<2\pi^2$.
We will
show (in Lemma \ref{any_monomial}) that all monomials ``$\z^\sk$"
can appear on the fibers of $E_m$, which is therefore holomorphically separable, and
``fiberwise convex" with respect to \inred{plurisubharmonic (}psh\inred{)} functions (Lemma \ref{fiberwise_convex}).
It also follows that there exists a continuous strictly psh
function $\psi$ on $E_m$. Then we show that $E_m$ is convex with respect
to continuous psh functions (Lemma \ref{psh_convex}).
From this and the existence of $\psi$ it follows by \cite{Nara}
({\em cf.} (0.2) in \cite{Siu}) that
if $\rho(M)=1$, or $\rho(M)>1$ and $m<\frac{2\pi^2}{\log \rho(M)}$, then $E_m$ is Stein.
For $\rho(M)=1$ we show the Steinness also when ``$m=2\infty$", i.e., the base is $\Cs$.

The case of $\rho(M)>1$ and $m= \tilde{m} = \frac{2\pi^2}{\log \rho(M)}$ is then deduced as follows:
From the above we obtain in particular an increasing sequence $(E^\nu)_{\nu\in\N}$
of Stein open subsets whose union is $E_{\tilde{m}}$ (e.g., $E^\nu = E_{\tilde{m}-1/\nu}$),
and such that for all $\nu$, $E^\nu$ is $\O(E^{\nu+1})$-convex.
This implies that $E_{\tilde{m}}$ is Stein ({\em cf.} (1.2) in \cite{Siu}).

\begin{lem}
\label{any_monomial}
Assume $m\log \rho(M)<2\pi^2$. Let $F=q\big(\{\wt\} \times \Csd\big)$ be any fiber of $E_m$.
For all \emph{$\k \in \Z^d$},
there exists a function $f$ on $E_m$
such that for all \emph{$\z \in \Csd,\quad f[\wt,\z]=\z^\sk$}.
\end{lem}

\begin{proof}
Fix any $\k\in \Z^d$.
Denote $\rho= \rho(M)$. We first assume that $\rho>1$.
Pick $\varepsilon>0$ so small that
\begin{equation}
\label{epsilon_is_small}
\big(\log (\rho + \varepsilon)\big) - \frac{2\pi^2}{m} < 0.
\end{equation}
Recall that $j.\k=M^j \k$.
By a Jordan form argument, there exists $C>0$ such that
$$\forall j \in \Z,\quad \norme{j.\k}\leq C (\rho +\varepsilon)^{|j|}.$$
Therefore
for all $j$ and $\z$, 
\begin{equation}
|\z^{j.\sk}| \stackrel{  {\scriptscriptstyle\text{by (\ref{ln_zk})} }}{=}
e^{\innerp{\mathbf{log}\ |\sz|}{j.\sk}}
\leq e^{ \alpha (\rho +\varepsilon)^{|j|}}
\leq e^{ \alpha e^{|j| \log(\rho +\varepsilon) }}
\label{growth_of_zjk}
\end{equation}
for some $\alpha>0$ that depends on $\z$ but not on $j$.

\medskip
Define $\Omega$ and $\Delta$ in $\O(S_m)$ by
$$\Omega(w)=e^{-2\cosh(\frac{2\pi^2}{m} w)} \text{ and } \Delta(w)=\frac{\sin \pi (w-\wt) }{\pi (w-\wt)}.$$
Key properties of these functions are:
\begin{equation}
 |\Omega(w)|\sim e^{-e^{\frac{2\pi^2}{m} |\mathit{Re}\ w|}}
 \text{ when $w\in S_m$ and } |\mathit{Re}\ w|\rightarrow \infty,
 \label{omega_s_growth}
\end{equation}
\begin{equation}
\text{for all } w\in S_m, \quad \{ \Delta(w+t)\ \vert\ t\in\R \} \text{ is bounded,}
\label{delta_bounded}
\end{equation}
\begin{equation}
\text{for all }j\in\Z, \quad \Delta(\wt + j)=
\begin{cases}
1 & \text{if $j=0$,}\\
0 & \text{if $j\neq 0$.}
\end{cases}
\label{delta}
\end{equation}
Property (\ref{omega_s_growth}) follows from the 
\modifinred{growth} of the function
$\cosh$ on any horizontal line contained
in the strip $S=\setofst{w\in\C}{\mathit{Im}\ |w|<\frac{\pi}{2}}$,
and that if $w\in S_m$ then 
$\frac{2\pi^2}{m}w \in S$.
Property (\ref{delta_bounded}) follows from the periodicity of $\sin$.

\smallskip
Now define $f\in \O(E_m) = \O^\Z\big(S_m \times \Csd\big)$ by
$$f(w,\z)= \frac{1}{\Omega(\wt)}  \sum_{j\in\Z} \Omega(w+j)\; \Delta(w+j)\; \z^{j.\sk}.$$
{\em Let's check that $f$ is well-defined:}
The $\Z$-invariance of the right hand side follows immediately from (\ref{jzk_equals_zjk}).
Fix $(w,\z)$. By (\ref{omega_s_growth}), for some $C_1>0$, for all $j$,
$$|\Omega(w+j)|\leq C_1\, e^{-e^{\frac{2\pi^2}{m} |\mathit{Re}\ w+j|}}.$$
As $|\mathit{Re}\ w+j|\geq |j|-|\mathit{Re}\ w|$,
by denoting $\beta = e^{-\frac{2\pi^2}{m} |\mathit{Re}\ w| } >0$ we obtain
$$|\Omega(w+j)|\leq C_1\, e^{- \beta e^{\frac{2\pi^2}{m} |j|}}.$$
From this, (\ref{growth_of_zjk}) and (\ref{delta_bounded}),
we obtain that for some $C_2>0$, for all $j$,
$$|\Omega(w+j) \; \Delta(w+j) \; \z^{j.\sk} |\leq
C_2\, e^{- \beta e^{\frac{2\pi^2}{m} |j|}}
e^{ \alpha e^{|j| \log(\rho +\varepsilon) }}.
%
$$
The right hand side of the above inequality behaves like the
{\em exponential} of
$$- \beta e^{\frac{2\pi^2}{m} |j|}
\underbrace{
\Big(
      1 - \frac{\alpha}{\beta} e^{|j| \big( \log(\rho +\varepsilon) -\frac{2\pi^2}{m} \big)}
\Big)
            }_{=x_j}
$$
It follows from (\ref{epsilon_is_small}) that $x_j \xrightarrow{|j|\rightarrow \infty} 1$.
Thus the series is pointwise absolutely converging because its general term decays
at a doubly exponential rate.

Moreover,
as (\ref{omega_s_growth}) and (\ref{delta_bounded}) hold locally uniformly
with respect to $w$, so does the convergence of the series. Thus $f$ is
holomorphic in $w$. Besides, for a fixed $w$, the series is
a (lacunary) Laurent series, so $f$ is holomorphic in $\z$.
By Hartogs's theorem, $f$ is a holomorphic function on $E_m$.

It follows from (\ref{delta}) that
for any $\z \in \Csd,\ f(\wt, \z)= \z^{\sk}$.

\medskip
Finally, if $\rho=1$, the growth of $\norme{j.\k}$ is now \modifinred{at most} 
polynomial of some degree $p'$.
Take an integer $p$ such that $2p>p'$.
We can apply a similar reasoning as above, with $\Omega$ defined by
$$\Omega(w)=e^{-w^{2p}},$$
which simultaneously fits the bill for any finite or infinite $m$.
\end{proof}

\begin{crl} Assume $m\log \rho(M)<2\pi^2$.
\label{fiberwise_convex}
Then
\begin{enumerate}
\item $E_m$ is holomorphically separable,
\item there exists a continuous strictly
psh function on $E_m$,
\item for any fiber $F$ there exists a continuous psh
function $\varphi_F$ on $E_m$ that restricts to an exhaustion on $F$.
\end{enumerate}
\end{crl}
\begin{proof}
The pull-back
to $E_m$ of $\iota: w\mapsto w$ on the base $A_m$ separates any two points that
do not lie on the same fiber.
Pick a fiber $F$. By Lemma \ref{any_monomial} applied to $F$ with
$\k_1=(1,0,\dots,0),\dots,\k_d=(0,\dots,0,1)$, we get functions
$g_1,\dots, g_{d}$; and with $-\k_1, \dots, -\k_d$ we get $g_{d+1},\dots,g_{2d}$.
The corresponding map $G: E_m \rightarrow \C^{2d}$
restricts on $F$ to a proper embedding because there is an isomorphism $H:\Csd\rightarrow F$
such that $GH(z_1,...,z_d)=(z_1,...,z_d,z_1\inv,...,z_d\inv)$.
This shows in particular that functions on $E_m$ separate
points of $F$.

Now let $\varphi_0$ be any continuous psh exhaustion on $\C^{2d}$.
Then $\varphi_F=\varphi_0 G$ is a continuous psh function on $E_m$ that restricts to an exhaustion
on $F$.

By the inverse mapping theorem,
the functions $\iota, f_1, \dots, f_d$ give local coordinates on a neighborhood
of any point of $F$. As $F$ was arbitrary, functions on $E_m$ give local
coordinates around every point, so one can construct a continuous strictly psh
function on $E_m$ as in \cite{Siu} (2.3) or \cite{Dem} I (6.17).
\end{proof}

From now on, for the case of $\rho(M)>1$, we look at all bundles $E_m$ with
$m<\tilde{m}=\frac{2\pi^2}{\log \rho(M)}$
as a family of open subsets of $E_{\tilde{m}}$.

Let $x=[w_1,\z_1]\in E_m$.
Notice that
if $v\in\C$ is so small that
$w_1+v$ belongs to $S_m$,
then for any $(w_2,\z_2)$ such that $[w_2,\z_2]=[w_1,\z_1]$,
($w_2+v$ also belongs to $S_m$ and) $[w_1 +v,\z_1]=[w_2+v,\z_2]$.
For convenience, we state and prove the following lemma, which is a mere reformulation
of an argument from \cite{Siu}.
\begin{lem}
\label{psh_convex}
Assume $m\log \rho(M)<2\pi^2$. Then $E_m$ is convex with respect
to continuous psh functions, i.e., for any closed discrete sequence
$\big(x_n=[w_n,\z_n]\big)_{n\in\N}$ there exists a continuous psh function $\varphi$ that
is unbounded on this sequence.
\end{lem}
\begin{proof}
If $p(x_n)$ has no accumulation point in the annulus $A_m$, then there
is a function $A_m$ whose pull-back on $E_m$ gives the desired $\varphi$.
So up to taking a subsequence, we can assume
that $p(x_n)$ converges to $w\in A_m$, and thus $(\z_n)_{n\in \N}$
can not have any accumulation point
in $D=\Csd$.

If $\rho(M)>1$, take $m'$ such that
$m < m'<\frac{2\pi^2}{\log \rho(M)}$, and
$\varepsilon>0$ such that: If $w\in S_m$ and $|v| \leq \varepsilon$, then $w+v\in S_{m'}$.
If $\rho(M)=1$, we assume that $m$ is infinite (i.e., $S_m=\C$)
and take any $\varepsilon>0$.

Denote $F=p\inv(w)$. Let $\varphi_F$ be the continuous psh function on $E_{m'}$
obtained from Corollary \ref{fiberwise_convex} (3).
Define a continuous psh function $\varphi$ on $E_m$ by
$$\varphi[w,\z]= \sup_{|v| \leq \varepsilon} \varphi_F[w+v,\z].$$
For $n$ big enough, $|w-w_n|<\varepsilon$, so $\varphi(x_n)\geq \varphi_F[w,\z_n]$.
Therefore $\varphi$ satisfies to the required properties.
\end{proof}
\inred{The Main Theorem is now proved.}

\subsection{Proof of Theorem \ref{fiber_Csd}.}
\label{proof_of_fiber_Csd}
The proof is simply a refinement of that of Sect. \ref{only_if}, based
on the extra assumptions we made. We omit the details.

Let $f\in\O(E_m)$. Let $\k$ be any nonzero element of $\Z^d$.
Take $\lambda_+$ and $\lambda_-$ eigenvalues of $M$ with respectively
maximal and minimal modulus. Denote $\mu_+=\log |\lambda_+|$
and $\mu_-=-\log |\lambda_-|$.

Take $\vb \in \C^d$ an eigenvector (for the linear action of
$M$ on $\C^d$) associated to $\lambda_+$.
We choose $\vb \in \R^d$ in case $\lambda_+$ is real. Define
$W=\mathit{Span}_{\R}\ \{\mathit{Re}\ \vb, \mathit{Im}\ \vb\}$.
Then $W$ is an $M$-invariant subspace of $\R^d$ on which $M$ acts
by a rotation-dilation of factor $\lambda_+$ (if $\lambda_+ \in \R$
then $\dim W=1$ and we have a dilation of factor $\lambda_+$).
As the characteristic polynomial of $M$ is irreducible, there
exists $\ub\in W$ such that $\innerp{\ub}{\k}>0$
({\em cf.} \cite{Fann-Wol} Lemma 3 (f)).

Take $\z_+ \in \Csd$ such that $\mathbf{log}\ |\z_+| =\ub$.
Now by a similar (but simpler) argument than for Lemma \ref{zd_and_kd},
we get an infinite $J_+\subset \N$ with bounded gaps such that
$|\z_+^{j.\sk}| > e^{\delta e^{j\mu_+}}$ for all $j\in J_+$.

Similarly, there exists $\z_- \in \Csd$ and an infinite
$J_-\subset -\N$ with bounded gaps such that
$|\z_-^{j.\sk}| > e^{\delta e^{|j|\mu_-}}$ for all $j\in J_-$.

\enlargethispage{1cm} Let $g_\sk$ be the coefficient of $f$ in the
Hartogs-Laurent series (\ref{grouped_coeffs}). Then, as in Sect.
\ref{only_if}, we obtain that $g_\sk$ vanishes identically. We
conclude that
$f(w,\z)=g_{\mathbf{0}}(w)$.\hfill$\square$\\

\section{Open problems.}

\cmt{Open problems section edited}
\subsection{Serre problem for bounded domains.}
Before C\oe ur\'e and L\oe b's counterexample, it was conjectured
by Siu that any bounded domain belonged to $\S$ ({\em cf.}
\ref{case_bounded_D}).
So far, the only known counterexamples
are (equivariant subsets of) Reinhardt domains
that generalize C\oe ur\'e and L\oe b's. This raises the problems,
already open in dimension two,
of the existence of other counterexamples, and of giving a characterization
of all bounded domains not in $\S$ (\cf \cite{CZ}).

\subsection{Siu's conjecture.}
For any of the known counterexamples to that conjecture,
our results say that, for a given
transition automorphism, a bundle with fiber $D$ over a
sufficiently thin annulus is Stein.
Note also that for another bundle, with fiber $\C^2$, an analogous result
is proved in \cite{Dem1}.
Therefore Siu's
conjecture is slightly rekindled, and this begs the question: Does
this interplay between fiber and base, that gives a Stein total
space, also happen for another ``sufficiently thin'' Stein base?
For other domains not in $\S$ (if they exist)?

\subsection{Stein and paralellizable manifolds.}
It is conjectured that a Stein and parallelizable $n$-manifold
(i.e., with trivial tangent bundle) can be realized as a Riemann
domain over $\Cn$ ({\em cf.} the nice survey \cite{For}).

The bundles $E_m(D,M)$ are parallelizable,
and when $m$ is so big that $E=E_m(D,M)$ is not Stein,
then $E$ is not a Riemann domain ({\em cf.} \cite{Siu}). So by
proving that for small $m$, $E$ is a Riemann domain, one would
obtain a good case study for the conjecture:
a continuous family of parallelizable manifolds (with a description of functions on each of them)
that degenerate
from Stein and Riemann domains to non-Stein and non-Riemann
domains.

\medskip
{\small {\bf Acknowledgements.} I am indebted to Al Baernstein \inred{and
Leonid Kovalev} for showing me that a stronger version of
Corollary 6.4 in \cite{Zaf}
followed from a ``classical" Phragm\'en-Lindel\"of-type argument.
A variation of their statement yields an important ingredient in Sect.
\ref{only_if}.
I am very grateful to Karl Oeljeklaus for many helpful conversations,
and to Charles Li for his listening and help. Thanks also go to
Alan Huckleberry for asking me a question that is answered by
Theorem \ref{fiber_Csd}.}

\end{document}